\documentclass[11pt,reqno,oneside]{amsart}
 \usepackage{latexsym}
 \usepackage{amssymb}
 \usepackage{amsfonts}
\usepackage{graphics}
\usepackage{mathcomp}
 
\thanks{University of Li\`ege, Institute of mathematics, Grande Traverse, 12 - B37, B-4000
 Li\`ege, Belgium email :  F.Radoux@ulg.ac.be}
\author[Radoux]{F. Radoux}
\date{\today} 
\title[Non-uniqueness]{Non-uniqueness of the natural and projectively equivariant quantization}

\newtheorem{lem}{Lemma}
\newtheorem{thm}[lem]{Theorem}
\newtheorem{prop}[lem]{Proposition}

\theoremstyle{remark}

\theoremstyle{definition}
\newtheorem{defi}{Definition}

\newcommand{\R}{\mathbb{R}}

\newcommand{\g}{\mathfrak{g}}

\begin{document}
\begin{abstract}
In \cite{DO}, the authors showed the existence and the uniqueness of a
$sl(m+1,\R)$-equivariant quantization in the non-critical situations. The
curved generalization of the $sl(m+1,\R)$-equivariant quantization is the
natural and projectively equivariant quantization. In \cite{Bord} and
\cite{moi}, the existence of such a quantization was proved in two different
ways. In this paper, we show that this quantization is not unique.

\end{abstract}
\maketitle
\noindent{\bf{Mathematics Subject Classification (2000) :}}  53B05, 53B10, 53D50, 53C10.\\
{\bf{Key words}} : Projective Cartan connections, differential operators,
natural maps, quantization maps.
\section{Introduction}
A quantization can be defined as a linear bijection from the space $\mathcal{S}(M)$ of symmetric contravariant tensor fields
on a manifold $M$ (also called the
space of \emph{Symbols}) to the space $\mathcal{D}_{\frac{1}{2}}(M)$ of
differential operators acting between half-densities.

It is known that
there is no natural quantization procedure. In other words, the spaces of
symbols and of differential operators are not isomorphic as representations of
$\mathrm{Diff}(M)$.

The idea of equivariant quantization, introduced by P. Lecomte and V. Ovsienko
in \cite{LO} is to reduce the group of local diffeomorphisms in the following
way.

They considered the case of the projective group
$PGL(m+1,\R)$ acting locally on the manifold $M=\R^m$ by linear fractional
transformations. They showed that the spaces of symbols and of differential
operators are canonically isomorphic as representations of $PGL(m+1,\R)$ (or
its Lie algebra $sl(m+1,\R)$). In other words, they showed that there exists a
unique \it{projectively equivariant quantization}. \rm{In} \cite{DO}, the authors
generalized this result to the spaces $\mathcal{D}_{\lambda\mu}(\R^{m})$ of
differential operators acting between $\lambda$- and $\mu$-densities and to
their associated graded spaces $\mathcal{S}_{\delta}$. They showed the
existence and uniqueness of a projectively equivariant quantization, provided
the shift value $\delta=\mu-\lambda$ does not belong to a set of critical
values. 

The problem of the $sl(m+1,\R)$-equivariant quantization on $\R^{m}$ has a
counterpart on an arbitrary manifold $M$. In \cite{Leconj}, P. Lecomte conjectured the existence of a quantization 
procedure depending on a torsion-free connection, that would be
 natural (in all arguments) and that would remain invariant by a projective
 change of connection.

After the proof of the existence of such a \it{Natural and equivariant
  quantization} \rm{given} by M. Bordemann in \cite{Bord}, we analysed in \cite{moi} the problem of this existence using
 Cartan connections. After these works, the question of the uniqueness of this
  quantization was not yet approached. The uniqueness of the $sl(m+1,\R)$-equivariant
  quantization in the non-critical situations did not imply the uniqueness of
  the solution in the curved case. The aim of this paper is to show that this
  solution is not unique, even in the non-critical situations, using the theory of Cartan connections.

\section{Fundamental tools}

For the sake of completeness, we briefly recall in this section the main notions and results
of \cite{moi}. We refer the reader to this reference or to \cite{Koba} for additional information. Throughout this note, we denote by $M$ a smooth, Hausdorff and
second countable manifold of dimension $m$.
\subsection{Natural and projectively equivariant quantization}
Denote by $\mathcal{F}_{\lambda}(M)$ the space of smooth
sections of the vector bundle of $\lambda$-densities. 

We denote by $\mathcal{D}_{\lambda,\mu}(M)$ the space of
differential operators from $\mathcal{F}_{\lambda}(M)$ to
$\mathcal{F}_{\mu}(M)$ and by $\mathcal{D}^k_{\lambda,\mu}$ the
space of differential operators of order at most $k$. 
If $\delta = \mu
- \lambda$, the associated space of \emph{symbols} will be called
$\mathcal{S}^k_{\delta}(M)$ and $\sigma$ will represent the \emph{principal symbol operator} from
$\mathcal{D}^k_{\lambda,\mu}(M)$ to $\mathcal{S}^k_{\delta}(M)$.

In these conditions, a \emph{quantization} on $M$ is a linear bijection $Q_M$ from
 the space of symbols
$\mathcal{S}_{\delta}(M)$ to the space of differential operators
$\mathcal{D}_{\lambda,\mu}(M)$ such that
\[\sigma (Q_M(S)) = S,\quad\forall S\in \mathcal{S}^k_{\delta}(M),\; \forall k
\in\mathbb{N}. \]
A \emph{natural quantization} is a quantization which
depends on a torsion-free connection and commutes with
 the action of diffeomorphisms.

A quantization $Q_M$ is \emph{projectively equivariant} if one has
$Q_M(\nabla) = Q_M(\nabla')$ whenever $\nabla$ and $\nabla'$ are projectively
equivalent torsion-free linear connections on $M$.

\subsection{Projective structures and Cartan projective connections}
We consider the group $G=PGL(m+1,\R)$ acting on the projective space.
 We denote by $H$ its isotropy subgroup at the origin. The group $H$ is the semi-direct product $G_0 \rtimes G_1$, where
$G_0$ is isomorphic to $GL(m,\R)$ and $G_1$ is isomorphic to $\R^{m*}$. The
Lie algebra associated to $H$ is $\g_0\oplus\g_1$. 

We recall that $H$ can be seen as a subgroup of the group of 2-jets 
$G_m^2$.\\
A \emph{Projective
  structure on $M$} is then a reduction of the second order frame bundle $P^2M$
to the group $H$. \\
The following result (\cite[p. 147]{Koba}) is 
the starting point of our method :
\begin{prop}[Kobayashi-Nagano]
There is a natural one to one correspondence between the projective
 equivalence classes of torsion-free linear connections on $M$
and the projective structures on $M$.
\end{prop}
We now refer the reader to \cite{Koba} for the definition of a projective
Cartan connection. Recall that if $\omega$ is a Cartan connection defined on a $H$-principal bundle $P$, then its
curvature $\Omega$ is defined as usual by
 \begin{equation}\label{curv} 
\Omega = d\omega+\frac{1}{2}[\omega,\omega].
\end{equation}
We can define from $\Omega$ a function $\kappa\in
C^{\infty}(P,\g_{-1}^{*}\otimes\g_{-1}^{*}\otimes\g)$ by :
$$\kappa(u)(X,Y):=\Omega(u)(\omega^{-1}(X),\omega^{-1}(Y)).$$
The \emph{Normal} Cartan connection has the following property (see
\cite[p. 136]{Koba}):
$$\sum_{i}\kappa_{jil}^{i}=0\quad\forall j,\forall l.$$

Now, the following result (\cite[p. 135]{Koba}) gives the relationship between  projective
structures and Cartan connections :
\begin{prop}
 A unique normal
 Cartan projective connection is
 associated to every projective structure $P$. This association is natural.
\end{prop}
The connection associated to a projective structure $P$ is called the normal
projective connection of the projective structure.
\subsection{Lift of equivariant functions}\label{Lift}

If $(V,\rho)$ is a representation of $\mathit{GL}(m,\R)$, then we can define
from it a representation $(V,\rho')$ of $H$ by projection (see \cite{moi} section 3). If $P$ is  a projective structure on $M$, the natural projection $P^2M\to
P^1M$ induces a projection $p :P\to P^1M$ and we have a well-known result:
\begin{prop} If $(V,\rho)$ is a representation of $GL(m,\R)$, then the map
\[p^* : C^{\infty}(P^1M,V)\to C^{\infty}(P,V) : f \mapsto f\circ p\]
defines a bijection from $C^{\infty}(P^1M,V)_{GL(m,\R)}$ to
$C^{\infty}(P,V)_{H}$.
\end{prop}

Subsequently, we will use the representation $\rho'_*$ of the Lie algebra of
$H$ on $V$. If we recall that this algebra is isomorphic to
 $gl(m,\R)\oplus \R^{m*}$
then we have 
\begin{equation}\label{rho}\rho'_* (A, \xi) = \rho_*(A),\quad\forall A\in
  gl(m,\R), \xi\in \R^{m*}.\end{equation}
Recall that $f\in C^{\infty}(P,V)_{G_1}$ if and only if
\begin{equation}\label{Invalg}
L_{h^*}f(u)=0,\quad\forall h\in\R^{m*}\subset sl(m+1,\R),
\forall u\in P.
\end{equation}

Finally, we recall the definitions of two operators used subsequently :
\begin{defi}
Let $(V, \rho)$ be a representation of $H$. If $f\in C^{\infty}(P,V)$, then $\nabla_{s}^{{\omega}^k} f
\in C^{\infty}(P,S^k\R^{m*}\otimes V)$ is defined by
\[\nabla_{s}^{{\omega}^k} f(u)(X_1,\ldots,X_k) =\frac {1}{k!}\sum_{\nu} L_{\omega^{-1}(X_{\nu(k)})}\circ\ldots\circ L_{\omega^{-1}(X_{\nu(1)})}f(u).\]
\end{defi}

If $(e_1,\ldots,e_m)$ is the canonical basis of $\R^m$ and if
$(\epsilon^1,\ldots,\epsilon^m)$ is the dual basis corresponding in $\R^{m*}$,
the \emph{divergence operator} is defined then by :
\[Div^{\omega} : C^{\infty}(P,S^k_{\delta}(\R^m))
\to C^{\infty}(P,S^{k-1}_{\delta}(\R^m)) :
S\mapsto \sum_{j=1}^m \nabla^{\omega}_{e_j}S(\epsilon^{j}).\]

\section{Non-uniqueness of the natural and projectively equivariant quantization}

First, one makes the following remark :
\begin{prop}
A natural and projectively equivariant quantization $Q_{M}$ is not unique if
and only if there is a nonzero natural projectively equivariant application
acting between $\mathcal{S}_{\delta}^{k}(M)$ and
$\mathcal{S}_{\delta}^{k-l}(M)$ for one $k$ and for one $l>0$.
\end{prop}
\begin{proof}
A quantization being a bijection, the non-uniqueness of a natural projectively
equivariant quantization is equivalent to the existence of two natural
projectively equivariant quantizations $Q$ and $Q'$ and of a natural
projectively equivariant application $T$ from $\mathcal{S}_{\delta}(M)$ to
$\mathcal{S}_{\delta}(M)$ different from the identity such that $Q'=Q\circ
T$. There is at least one $k$ such that the restriction of $T$ to
$\mathcal{S}_{\delta}^{k}(M)$ is different from the identity. As a
quantization must preserve the principal symbol, the projection of this
restriction on $\mathcal{S}_{\delta}^{k}(M)$ must be equal to the
identity. The projections of the restriction on
$\mathcal{S}_{\delta}^{k+l}(M)$, with $l>0$, must be equal to zero and one can
conclude.  
\end{proof}

The construction of the applications discussed in the previous result is based
on the Weyl tensor. Let us first recall its definition.

\subsection{The Weyl tensor}
If we denote by $\omega$ the normal Cartan connection associated to a
projective structure $P$, the function $\kappa$ induced by its curvature has
an important property of invariance :
\begin{prop}
If $h\in H$, $u\in P$ and $X,Y\in\g_{-1}$, the function $\kappa\in
C^{\infty}(P,\g_{-1}^{*}\otimes\g_{-1}^{*}\otimes\g)$ satisfies :
\begin{equation}\label{kappa}
\kappa(X,Y)(uh)=Ad(h^{-1})\kappa(Ad(h)X,Ad(h)Y)(u).
\end{equation}
The function $\kappa_{0}$ is accordingly $H$-equivariant. It represents
then a tensor of type $\left(\begin{array}{c}1\\3\end{array}\right)$ on $M$
  that is called the \it{Weyl tensor}.
\end{prop}
\begin{proof}
The relation (4) is proved in \cite{Slo1} page 44.

If one considers the components according to $\g_{0}$ of the two hands of
(\ref{kappa}), one has :
$$\kappa_{0}(X,Y)(uh)=\rho^{\R^{m}\otimes\R^{m*}}(h^{-1})\kappa_{0}(\rho^{\R^{m}}(h)X,\rho^{\R^{m}}(h)Y)(u),$$
where $\rho^{\R^{m}\otimes\R^{m*}}$ and $\rho^{\R^{m}}$ denote respectively
the actions of $H$ on $\R^{m}\otimes\R^{m*}$ and $\R^{m}$. Indeed, the
components according to $\g_{-1}$ of $Ad(h)X$ and $Ad(h)Y$ are respectively 
$\rho^{\R^{m}}(h)X$ and $\rho^{\R^{m}}(h)Y$. Moreover, the fact that
$\kappa_{-1}=0$ implies that the component according to $\g_{0}$ of
$$Ad(h^{-1})\kappa(\rho^{\R^{m}}(h)X,\rho^{\R^{m}}(h)Y)(u)$$
is equal to
$$\rho^{\R^{m}\otimes\R^{m*}}(h^{-1})\kappa_{0}(\rho^{\R^{m}}(h)X,\rho^{\R^{m}}(h)Y)(u).$$ 
\end{proof}

\subsection[Applications naturelles projectivement invariantes]{Construction
  of natural and projectively equivariant applications}

\rm{If} $j\geq 2$ and if $\sigma$ is a permutation of $\{1,\ldots,j\}$ such
that \\
$\sigma(l)\neq l\quad\forall l\in\{1,\ldots,j\}$, we define a $H$-equivariant function $W\in
C^{\infty}(P,S^{2j}\R^{m*})$ in the following way :
$$W(u)(e_{i_1},\ldots,e_{i_{2j}}):=\sum_{\nu}\sum_{r_{1},\ldots,r_{j}}{\kappa_{0}(u)}_{i_{\nu(1)}i_{\nu(2)}r_{\sigma(1)}}^{r_1}\ldots{\kappa_{0}(u)}_{i_{\nu(2j-1)}i_{\nu(2j)}r_{\sigma(j)}}^{r_j}.$$

The following lemma allows to calculate the failure of equivariance of the
iterated invariant differentials of the function $W$ :
\begin{lem}
One has the following formula :
$$L_{h^{*}}\nabla_{s}^{\omega^k}W=-k(k+4j-1)h\vee(\nabla_{s}^{\omega^{k-1}}W),$$
for all $h\in\g_{1}$.
\end{lem}
\begin{proof}
The proof is similar to the proof of the proposition 10 of \cite{moi}. If $k=0$,
the formula is true. One proceeds then by induction. If $X\in\R^{m}$, 
$$(L_{h^{*}}\nabla_{s}^{\omega^k}W)(X,\ldots,X)$$
is equal to
$$L_{\omega^{-1}(X)}L_{h^{*}}(\nabla_{s}^{\omega^{k-1}}W)(X,\ldots,X)+(L_{[h,X]^{*}}(\nabla_{s}^{\omega^{k-1}}W))(X,\ldots,X).$$
By induction, the first term is equal to
$$-(k-1)(k+4j-2)(h\vee(\nabla_{s}^{\omega^{k-1}}W))(X,\ldots,X).$$
Concerning the second term, one obtains 
\[(\rho_{*}((h\otimes X) +\langle h, X\rangle Id)(\nabla_s^{{\omega}^{k-1}}W))(X,\ldots,X).\]
The result comes then from the definition of $\rho_*$.
\end{proof}
One can then construct natural projectively invariant applications between
spaces of symbols :
\begin{thm}\label{princ} 
If $S\in\mathcal{S}_{\delta}^{k}(M)$ and $l\geq 2j$, the multiples of the application
\begin{equation*}S\mapsto p^{*^{-1}}(\sum_{r=0}^{l-2j} C_{k,l,r} \langle Div^{\omega^r}
p^*S,\nabla_s^{\omega^{l-r-2j}}W\rangle)\end{equation*}
are natural and projectively equivariant if
\[C_{k,l,r} =\frac{(l+2j-1)!}{(m+1)^{r}(l+2j-1-r)!\gamma_{2k-1}\cdots
  \gamma_{2k-r}}\left(\begin{array}{c}l-2j\\r\end{array}\right),\forall r\geq 1,\;C_{k,l,0}=1.\]
\end{thm}
\begin{proof}
The proof is similar to the proof of the theorem 11 of \cite{moi}. Thanks to the
proposition 4 and to the lemma 7 of \cite{moi}, it suffices to check that the function
\begin{equation*}\sum_{r=0}^{l-2j} C_{k,l,r} \langle Div^{\omega^r}
p^*S,\nabla_s^{\omega^{l-r-2j}}W\rangle\end{equation*}
is $\g_{1}$-equivariant. This is true thanks to the previous lemma, to the
proposition 9 of \cite{moi} and to the fact that the following relation is
satisfied
:
\begin{equation*}C_{k,l,r}r(m+2k-r-(m+1)\delta)=C_{k,l,r-1}(l-r-2j+1)(l-r+2j).\end{equation*}
The application given in the theorem is projectively equivariant by definition
of $\omega$. It is natural too : it follows from the naturality of the
association of a projective structure $P\to M$ endowed with a normal Cartan
connection $\omega$ to a class of projectively equivalent torsion-free
connections on $M$.
\end{proof}
\noindent{\bf{Remarks :}}
\begin{itemize}
\item The applications that we have given are some examples of natural projectively
equivariant applications between spaces of symbols. A complete description of
the set of these applications seems however rather difficult.  
\item One can show ``by hand'' that the natural and projectively
  equivariant quantization is unique up to the third order in the non-critical
  situations. It suffices to consider all the natural applications between
  $\mathcal{S}^{k}(M)$ and $\mathcal{S}^{k-l}(M)$ (with $1\leq l\leq 3$) and to
  show that there is not linear combination of these maps that is projectively
  equivariant in the non-critical situations.
\item Using methods described in \cite{Slo1} and \cite{formules}, one could derive explicit formulas
  for the applications of the theorem \ref{princ}. At the fourth and fifth
  orders, if we denote by $T$ the equivariant function on $P^{1}M$ corresponding to
  $W$ with $j=2$, the applications of the theorem \ref{princ} are equal
  respectively to  
$$\langle S,T\rangle$$
and to 
$$\langle S,\nabla_s T\rangle+\frac {8}{(m+1)\gamma_{2k-1}}\langle Div S,T\rangle.$$
\end{itemize}

\section{Acknowledgements}

 I thank the Belgian FRIA for my Research Fellowship.

\bibliographystyle{plain} \bibliography{unic}

\begin{thebibliography}{1}

\bibitem{Bord}
M.~Bordemann.
\newblock Sur l'existence d'une prescription d'ordre naturelle projectivement
  invariante.
\newblock {\em Submitted for publication}, math.DG/0208171.

\bibitem{Slo1}
A.~{\v{C}}ap, J.~Slov{\'a}k, and V.~Sou{\v{c}}ek.
\newblock Invariant operators on manifolds with almost {H}ermitian symmetric
  structures. {I}. {I}nvariant differentiation.
\newblock {\em Acta Math. Univ. Comenian. (N.S.)}, 66(1):33--69, 1997.

\bibitem{DO}
C.~Duval and V.~Ovsienko.
\newblock Projectively equivariant quantization and symbol calculus:
  noncommutative hypergeometric functions.
\newblock {\em Lett. Math. Phys.}, 57(1):61--67, 2001.

\bibitem{Koba}
Shoshichi Kobayashi.
\newblock {\em Transformation groups in differential geometry}.
\newblock Springer-Verlag, New York, 1972.
\newblock Ergebnisse der Mathematik und ihrer Grenzgebiete, Band 70.

\bibitem{LO}
P.~B.~A. Lecomte and V.~Yu. Ovsienko.
\newblock Projectively equivariant symbol calculus.
\newblock {\em Lett. Math. Phys.}, 49(3):173--196, 1999.

\bibitem{Leconj}
Pierre B.~A. Lecomte.
\newblock Towards projectively equivariant quantization.
\newblock {\em Progr. Theoret. Phys. Suppl.}, (144):125--132, 2001.
\newblock Noncommutative geometry and string theory (Yokohama, 2001).

\bibitem{moi}
P.~Mathonet and F.~Radoux.
\newblock Natural and projectively equivariant quantizations by means of
  {C}artan connections.
\newblock {\em Lett. Math. Phys.}, 72(3):183--196, 2005.

\bibitem{formules}
F.~Radoux.
\newblock Explicit formula for the natural and projectively equivariant
  quantization.
\newblock {\em Submitted for publication, math.DG/0606522}.

\end{thebibliography}
\end{document}